\documentclass{amsart}
\usepackage{graphicx}
\newtheorem{thm}{Theorem}[section]
\newtheorem{cor}[thm]{Corollary}
\newtheorem{lem}[thm]{Lemma}

\theoremstyle{definition}
\newtheorem{defn}[thm]{Definition}
\theoremstyle{remark}
\newtheorem{rem}[thm]{Remark}
\newtheorem{question}[thm]{Question}
\numberwithin{equation}{section}

\newcommand{\Real}{\mathbb R}

\newcommand{\sinc}{\mathrm{sinc\:}}
\newcommand{\arccot}{\mathrm{arccot\:}}
\begin{document}

\title[]{The Amplitude Modulation transform}%
\author{Igor Rivin}%
\address{Mathematics Department, Temple University}%
\email{rivin@math.temple.edu}%

\thanks{The author was supported by a grant from the National Science Foundation.
He would like to thank Warren D. Smith for enlightening discussions}%

\begin{abstract}
Motivated by the study of the local extrema of $\sin(x)/x$ we
define the \emph{Amplitude Modulation} transform of functions
defined on (subsets of) the real line. We discuss certain
properties of this transform and invert it in some easy cases.
\end{abstract}
\maketitle
\section*{Introduction}

This note has been motivated by the following question:

\begin{quotation}
Let $0 = x_0 < x_1 < \cdots < x_n < \cdots$ be the sequence of
local maxima of the \emph{sinc} function $\sinc(x) = \sin(x)/x.$
Is the sequence $1 = \sinc(x_0), \sinc(x_1), \dots, \sinc(x_n),
\dots$ decreasing?
\end{quotation}

This question is not difficult to answer. Indeed, at a critical
point $x_i,$
$$\sinc^\prime(x_i) = \frac{x_i \cos(x_i) - \sin(x_i)}{x_i^2}=0,$$
which implies that
\begin{equation}
\label{valeq} \cos(x_i) = \frac{\sin(x_i)}{x_i} = \sinc(x_i).
\end{equation}
One can also write the above equation as:
\begin{equation}
\label{pteq}
x_i = \tan(x_i),
\end{equation}
or, equivalently:
\begin{equation}
\label{pteq2}
 \arctan(x_i) = x_i.
\end{equation}
Combining equations (\ref{valeq}) and (\ref{pteq2}), we obtain:
\begin{equation}
\label{basiceq}
 \sinc(x_i) = \cos(\arctan(x_i)) =
\frac{1}{\sqrt{1+x_i^2}},
\end{equation}
so the decrease of $\sinc(x_i)$ is immediate.
\section{The \emph{Amplitude Modulation} transform}
The formula (\ref{basiceq}) suggest the following:
\begin{defn}
The \emph{Amplitude Modulation} transform $\mathcal{AM}(f)$ of a
function $f:\Real\to\Real$ is the set of functions whose values at
the critical points of $f(x) \sin(x)$ agrees with those of $f(x)
\sin(x).$
\end{defn}

 \begin{rem}
 In fact, if in the definition of the $\mathcal{AM}$ transform we replace the
multiplier $\sin(x)$ by $\sin(x+k),$ we obtain the same function
$\mathcal{A}\mathcal{M}f.$ This is an observation of W.~D.~Smith,
and it allows us to replace the definition above with the more
pleasant definition below:
\end{rem}

\begin{defn}
The \emph{Amplitude Modulation} transform $\mathcal{AM}(f)$ of a
function $f:\Real\to\Real$ is the function whose values at the
critical points of $f(x) \sin(x+k)$ agrees with those of $f(x)
\sin(x+k)$ for \emph{all} values of the phase parameter $k.$
\end{defn}
 The discussion in the Introduction can thus be summarized as
follows:
\begin{thm}
The function $\frac{1}{\sqrt{1+x^2}}$ is the $\mathcal{AM}$
transform of $\frac{1}{x}.$
\end{thm}
To get an analogous result for a general function $f(x),$ we
perform the same sort of computation as in the Introduction:

(We will use the notation $f_x$ for the derivative of $f$ for
typographical reasons.) The critical points of $f(x) \sin(x)$ are
the points where:
\begin{equation*}
\frac{d f(x) \sin(x)}{d x} = 0.
\end{equation*}
Expanding, we see that $f(x) \cos(x) + f_x(x) \sin(x) = 0,$ and so
\begin{equation*}
\cot(x) = -\frac{f_x(x)}{f(x)},
\end{equation*}
so that
\begin{equation}
\label{eqval2}
 x = \arccot \left(-\frac{f_x(x)}{f(x)}\right),
\end{equation}
while
\begin{equation}
\label{eqval3}
 f(x) \sin(x) = - \tan(x) \sin(x) f_x(x).
\end{equation}
Combining Eq. (\ref{eqval2}) and Eq. (\ref{eqval3}) we see that at
the critical points:
\begin{equation}
\label{maineq}
f(x) \sin(x) = \pm \frac{f^2(x)}{\sqrt{f^2(x) + f_x^2(x)}},
\end{equation}
which we can summarize in
\begin{thm}[Theorem-Definition]
The function $\mathcal{AM}(f)$ is defined by
\begin{equation*}
\mathcal{AM}(f)(x) =  \frac{f^2(x)}{\sqrt{f^2(x) + f_x^2(x)}},
\end{equation*}

 \end{thm}

Here are some examples: As we have seen before, if $f(x) = 1/x,$
then
$$\mathcal{AM}(f)(x) = \frac{1}{\sqrt{1+ x^2}}.$$
\begin{figure}
\includegraphics[scale=1.0]{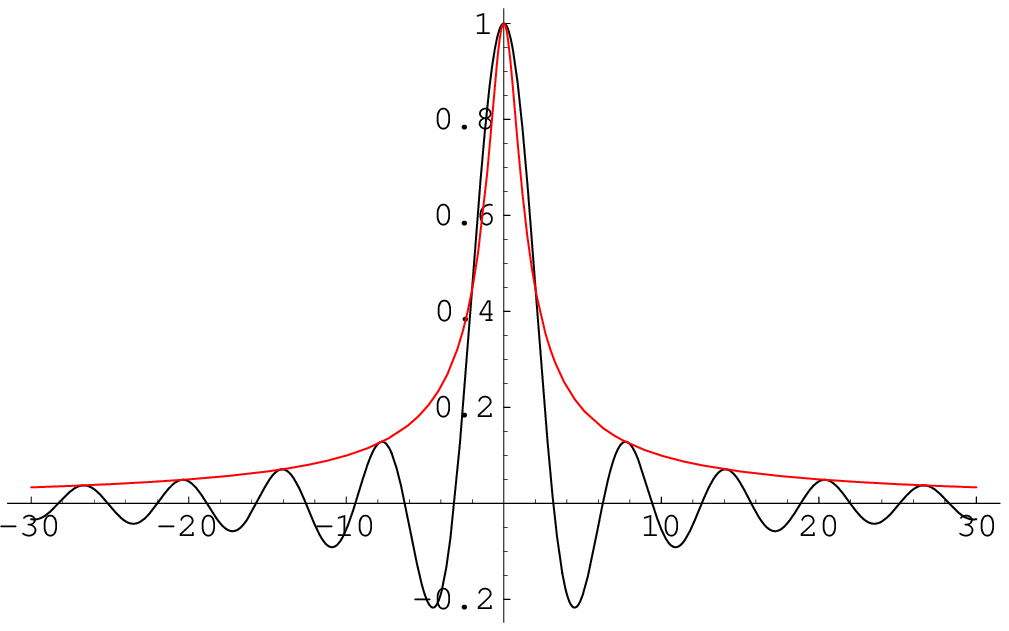}\\
\caption{$f(x) = 1/x.$}\label{sincfig}
\end{figure}

If $f(x) = x^\alpha,$ then
\begin{equation*}
\mathcal{AM}(f)(x)
  = \frac{x^{\alpha+1}}{\sqrt{x^2+\alpha^2}}.
\end{equation*}
\begin{figure}
  \includegraphics[scale=1.0]{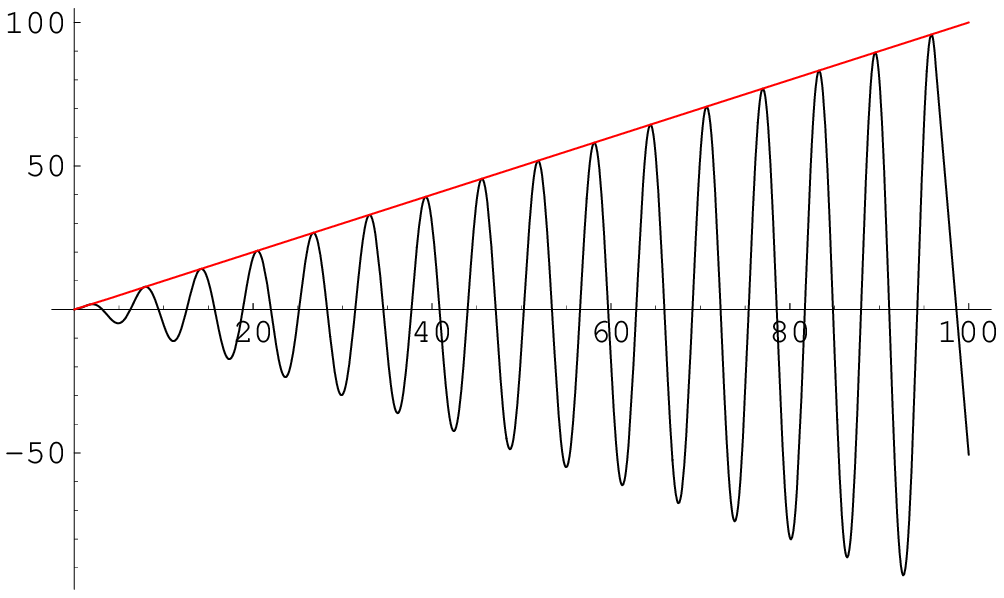}\\
  \caption{$f(x) = x.$}\label{exfig}
\end{figure}

If $f(x) = \exp(x),$ then
\begin{equation*}
\mathcal{AM}(f)(x) = \frac{\exp(x)}{\sqrt{2}},
\end{equation*}
while if $f(x) = \exp(-x),$ then
\begin{equation*}
\mathcal{AM}(f)(x) = \frac{\exp(-x)}{\sqrt{2}}.
\end{equation*}
\begin{figure}
  \includegraphics[scale=1.0]{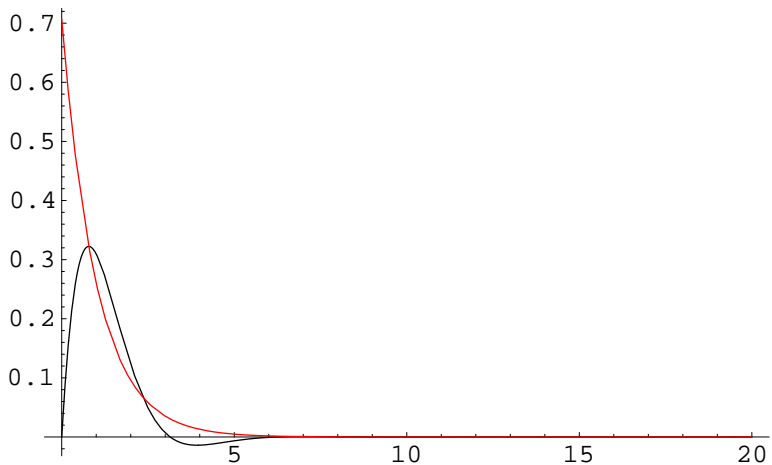}\\
  \caption{$f(x)=\exp(-x).$}\label{expfig}
\end{figure}

If $f(x) = \exp(g(x)),$ then
\begin{equation}\label{expform}
\mathcal{AM}(f)(x) = f(x) \frac{1}{\sqrt{1+g_x^2}}.
\end{equation}
 \section{Some algebraic observations}
 We will need to recall a definition:
 \begin{defn}
 A function $y = f(x)$ is called \emph{algebraic} if there exists
 a two-variable polynomial $P,$ such that $P(x, y) = 0.$
 \end{defn}
And some well-known results:
\begin{thm}
\label{elim} If $z$ is an algebraic function of $y$ and $y$ is an
algebraic function of $x,$ then $z$ is an algebraic function of
$x.$ Further, if $y_1$ and $y_2$ are algebraic functions of $x,$
then so are $y_1 y_2$ and $y_1 + y_2.$ Finally, if $y$ is an
algebraic function, then so is $y^\prime.$
\end{thm}
\begin{proof}[Proof Sketch and lightning introduction to
elimination theory]
The proofs of all the assertions follow from the following basic
fact: two univariate polynomials $P_1$ and $P_2$ over a domain $R$
with unity have no common zeros if and only if their greatest
common divisor is $1,$ or, equivalently, there exist polynomials
$Q_1$ and $Q_2,$ such that
\begin{equation}
\label{gcdeq} Q_1 P_1 + Q_2 P_2 = 1.
\end{equation}
Since Eq. (\ref{gcdeq}) is a system of linear equations for the
coefficients of $Q_1$ and $Q_2,$ the existence of $Q_1$ and $Q_2$
as above is easily seen to be equivalent to the non-vanishing of
the determinant of the linear system. This determinant is the
so-called \emph{resultant} of the polynomials $P_1$ and $P_2.$
Now, if we have two polynomial equations $P(x, y) = 0$ and $Q(x,
y) = 0,$ they can regarded as two polynomials in $y$ whose
coefficients are polynomials in $x,$ and so the set of
$x$-coordinates of the points in the common zero-set of $P$ and
$Q$ all have the property that at those points, the resultant of
the two equations vanishes. The resultant is a polynomial in $x,$
and so $y$ has been eliminated from consideration, hence the name
``elimination theory.'' Now, to proceed with the proof of the
Theorem \ref{elim}: If $z$ is an algebraic function of $y$ and $y$
is an algebraic function of $x,$ then there are equations $P(x,
y)=0$ and $Q(y, z)=0.$ Eliminating $y$ from the two equations, we
see that $z$ is algebraic. If $y_1$ and $y_2$ are algebraic, and
$y = y_1 y_2,$ then we have the the three equations (the ones
satisfied by $y_1$ and $y_2$ and $y = y_1 y_2.$ We can eliminate
first $y_1$ and then $y_2,$ to show that $y$ is algebraic,
similarly with $y_1 + y_2.$

Finally, to show that the derivative is algebraic, we
differentiate $P(x, y)=0$ implicitly, to obtain $Q(x, y, y^\prime)
= 0.$ Eliminating $y$ from the two equations we obtain the
algebricity of $y^\prime.$
\end{proof}
\begin{rem} For considerably more detail on the subject of
elimination, please see \cite{wang}. \end{rem}
 As corollaries
of the above Theorem, we see that
\begin{cor}
If $f(x)$ is an algebraic function, then so is $\mathcal{S}(f).$
If $f, g$ are algebraic functions, then if $h(x) = f(x)
\exp(g(x)),$ it follows that $\mathcal{AM}(h)(x) = k(x) h(x),$
where $k$ is algebraic.
\end{cor}

\section{Inverse problems}
The first obvious inverse problem is the following:
\begin{quotation}
Which functions $g(x)$ are $\mathcal{AM}$ transforms?
\end{quotation}
Construction of the inverse transform is equivalent to the
solution of the ODE
\begin{equation}
\label{invode}
 f^\prime(x) = \pm f(x)\sqrt{\frac{f^2(x)}{g^2(x)} -
1}.
\end{equation}
The choice of plus or minus is already troubling, as is the fact
that the right hand side is frequently not Lipschitz, so the usual
Picard existence theorem for ODE does not apply everywhere, and
uniqueness fails spectacularly: the functions $1/\sin x$ and $1$
have the same $\mathcal{AM}$ transform. This example also
demonstrates that the initial value problem can develop
singularities in finite time. Nevertheless, some things can be
said. First:

\begin{lem}
\label{invlem}
 The transformed function $\mathcal{AM}(f)$ has a
critical point whenever $f$ has a critical point. Furthermore, at
such a critical point $x,$ $\mathcal{AM}(f)(x) = f(x),$ and the
last equality only holds at a critical point of $f.$
\end{lem}

\begin{proof}
A simple computation.
\end{proof}

We can thus simplify our life by attempting to solve \ref{invode}
on an interval $[a, b]$ where $g$ is monotone, and in addition, $0
< g(x) < M.$ We can pick between the two equations:
\begin{eqnarray}
\label{invodep} f^\prime(x) &=& f(x) \sqrt{\frac{f^2(x)}{g^2(x)} -
1}\\
\label{invodem} f^\prime(x) &=& - f(x) \sqrt{\frac{f^2(x)}{g^2(x)}
- 1}
\end{eqnarray}

Now, by the Lemma \ref{invlem} we know that we have local
existence and uniqueness of solutions, and so the only thing we
need check is that singularities do not develop in finite time. To
do this we analyze two separate cases:
\begin{itemize}
\item{Case 1.} $g(x)$ is \emph{decreasing} on $[a, b].$ In this
case we take Eq.~(\ref{invodem}). Local existence and uniqueness
is assured by the Picard theorem (see \cite[Chapter 1]{taylor1}).
We pick the initial value $f(a)$ at will (as long as it is bigger
than $g(a).$)  Since $f^\prime(x)$ is always negative we know that
$f(x) < f(a),$ and since we know that $g(x) > 0$ on $[a, b]$ we
know that $f(x) > g(x) > 0,$ so it follows that we have a solution
on $[a, b].$
\item{Case 2.} $g(x)$ is \emph{increasing}. In this case we start
at the right endpoint $b,$ and use Eq.~(\ref{invodep}), and then
construct the solution going right to left. The reasoning in Case
1 goes through verbatim.
\end{itemize}

What happens if $g(x)$ has critical points on $[a, b]$? In that
case, it is fairly obvious that we can construct a ``weak
inverse,'' but anything more seems to require much more work.
In case the reader is dissatisfied with the nonexplicit nature of
our construction, (s)he will perhaps be mollified by the
observation that the following problem can be solved explicitly:

\begin{quotation}
Given a function $r(x) > 1$ on $[a, b],$ construct a positive
$f(x),$ such that $f(x)/\mathcal{AM}(f)(x) = r(x).$
\end{quotation}
Using Eq.~(\ref{expform}), it is easy to see that
\begin{equation}
f(x) = \exp\left(\pm C \int_a^x dt \sqrt{r^2(t) - 1} \right)
\end{equation}
is the desired solution.
\section{Questions}
The most natural question is:
\begin{question}
Given a smooth $g(x),$ is there a natural way to construct an
$f(x)$ such that $g(x) = \mathcal{AM}(f)$
\end{question}
Changing categories:
\begin{question}
Is $\mathcal{AM}$ invertible on the set of algebraic functions?
\end{question}
or
\begin{question}
Suppose $f$ satisfies a first order linear differential equation
with algebraic coefficients. Is the same always true of
$\mathcal{AM}(f)?$
\end{question}

\bibliographystyle{amsplain}

\end{document}